\documentclass{amsart}
\usepackage{amssymb}
\usepackage{amsmath}
\usepackage{amscd}
\usepackage[left=3cm, right=2.5cm, top=2cm, bottom=2cm]{geometry}
\usepackage[dvipsnames]{xcolor}
\usepackage{mathtools}

\newtheorem{prop}{Proposition}[section]

\newtheorem{lemma}[prop]{Lemma}
\newtheorem{theorem}[prop]{Theorem}
\newtheorem{cor}[prop]{Corollary}

\theoremstyle{definition}

\newcommand{\NN}{{\mathbb N}}

\DeclareFontEncoding{OT2}{}{} % to enable usage of cyrillic fonts

\begin{document}

\title{Some New Congruences for $\ell$-Regular Multipartitions}
\author{Srilakshmi Krishnamoorthy}
\author{Abinash Sarma}

\address{Indian Institute of Science Education and Research Thiruvananthapuram, Maruthamala P.O., Vithura, Thiruvananthapuram-695551, Kerala, India.}
\email{srilakshmi@iisertvm.ac.in}
\address{Indian Institute of Science Education and Research Thiruvananthapuram, Maruthamala P.O., Vithura, Thiruvananthapuram-695551, Kerala, India.}
\email{sarmaabinash15@iisertvm.ac.in}

\begin{abstract}
For a positive integer $n$, let $B_{\ell_1,\dots,\ell_r}(n)$ denote the number of $(\ell_1,\ell_2,\cdots,\ell_r)$-regular multipartitions of $n$. If $\ell_1=\ell_2=\cdots=\ell_r=\ell$, then we denote $B_{\ell_1,\dots,\ell_r}(n)$ as $B_\ell^{(r)}(n)$. In this paper, we prove several infinite families of congruences satisfied by $B_\ell^{(r)}(n)$ for different values of $\ell$ and $r$.
\end{abstract}

\maketitle

\section{\Large{\textbf{Introduction}}}

Given a positive integer $n$, a partition of $n$ is said to be $\ell$-regular, for a fixed positive integer $\ell$, if none of its parts is divisible by $\ell$. Let $n=a_1+a_2+\cdots+a_r$ be a partition of $n$. Suppose, in turn, $\lambda^{(i)}$ is a partition of $a_i$ for each $i\in\{1,2,\dots,r\}$. We call the ordered tuple $(\lambda^{(1)},\lambda^{(2)},\dots,\lambda^{(r)})$ as a $r$-multipartition of $n$; if the value of $r$ is obvious from the context, we simply call it a multipartition of $n$. In particular, a $2$-multipartition of $n$ is called a bipartition of $n$. Let us fix $r$ positive integers, namely, $\ell_1,\ell_2,\dots,\ell_r$ such that $\lambda^{(i)}$ is a $\ell_i$-regular partition. Then $(\lambda^{(1)},\lambda^{(2)},\dots,\lambda^{(r)})$ is said to be a $(\ell_1,\ell_2,\dots,\ell_r)$-regular multipartition of $n$.

Let $B_{\ell_1,\dots,\ell_r}(n)$ denote the number of $(\ell_1,\ell_2,\dots,\ell_r)$-regular multipartitions of $n$. By convention, we assume that $B_{\ell_1,\dots,\ell_r}(0)=1$ for any $r$ and any $(\ell_1,\ell_2,\dots,\ell_r)$. Note that the generating function for $B_{\ell_1,\dots,\ell_r}(n)$ is given by
\begin{align} \label{eq:1}
    \sum_{n=0}^\infty B_{\ell_1,\dots,\ell_r}(n)q^n=\frac{E_{\ell_1}E_{\ell_2}\cdots E_{\ell_r}}{E_1^r},
\end{align}
where $E_k:=\prod_{m=1}^\infty(1-q^{km})$ for $k\in\mathbb N$. Note that, for any prime number $p$ and any $k\in\NN$, we have
\begin{align} \label{eq:3}
    E_{kp}\equiv E_k^p\pmod{p}.
\end{align}
Suppose $\ell_1=\ell_2=\dots=\ell_r=\ell$. Then we say that $(\lambda^{(1)},\lambda^{(2)},\dots,\lambda^{(r)})$ is a $\ell$-regular $r$-multipartition of $n$. In this case, we simplify the notation for $B_{\ell_1,\dots,\ell_r}(n)$ as $B_\ell^{(r)}(n)$. Now from (\ref{eq:1}), we have
\begin{align} \label{eq:2}
    \sum_{n=0}^\infty B_\ell^{(r)}(n)q^n=\frac{E_\ell^r}{E_1^r}.
\end{align}

Recently, for different values of $\ell$, many properties of $\ell$-regular partition function have been discussed by several mathematicians in \cite{adiga&dasappa,ahmed&baruah,carlson&webb,cui&gu,dao&liu&yan,furcy&penniston,hou&sun&zhang,ranganatha,webb,xia}. For recent results for $(\ell_1,\ell_2)$-regular bipartition function, the reader may refer to \cite{adiga&ranganatha,dou,kathiravan,kathiravan&fathima,lin,lin2,wang}. Saikia and Boruah \cite{saikia&boruah} proved certain congruences for $B_\ell^{(3)}(n)$ taking $\ell\in\{2,3,4,5\}$ and Gireesh and Naika \cite{gireesh&naika} proved for $B_3^{(3)}(n)$. Chern, Tang, and Xia \cite{chern&tang&xia} established infinite families of congruences for $B_7^{(3)}(n)$ modulo powers of $7$. Baruah and Das \cite{baruah&das} established the generating functions for $B_9^{(3)}(n)$ and $B_{27}^{(3)}(n)$. Very recently, Murugan and Fathima \cite{murugan&fathima} studied congruences for $B_3^{(6)}(n)$.

The following are the main results of this paper.

\begin{theorem} \label{thm1}
    For any $n\geq0$ and $t\geq0$, we have
    \begin{itemize}
        \item[(i)] $B_3^{(12)}\left(3p_1^2\cdots p_{t+1}^2n+\left(p_1^2\cdots p_t^2p_{t+1}(p_{t+1}+3j)-1\right)\right)\equiv0\pmod3,$ where $p_i$ is a prime number such that $p_i\equiv2\pmod3$ for $i\in\{1,\cdots,t+1\}$ and $p_{t+1}\nmid j$.

        \item[(ii)] $B_3^{(15)}\left(3p_1^2\cdots p_{t+1}^2n+\frac{p_1^2\cdots p_t^2p_{t+1}(5p_{t+1}+12j)-5}{4}\right)\equiv0\pmod3,$ where $p_i\neq3$ is a prime number such that $p_i\equiv3\pmod4$ for $i\in\{1,\cdots,t+1\}$ and $p_{t+1}\nmid j$.

        \item[(iii)] $B_5^{(10)}\left(5p_1^2\cdots p_{t+1}^2n+\frac{5\big(p_1^2\cdots p_t^2p_{t+1}(p_{t+1}+3j)-1\big)}{3}\right)\equiv0\pmod5,$ where $p_i\neq5$ is a prime number such that $p_i\equiv2\pmod3$ for $i\in\{1,\cdots,t+1\}$ and $p_{t+1}\nmid j$.

        \item[(iv)] $B_7^{(7)}\left(7p_1^2\cdots p_{t+1}^2n+\frac{7\big(p_1^2\cdots p_t^2p_{t+1}(p_{t+1}+4j)-1\big)}{4}\right)\equiv0\pmod7,$ where $p_i\neq7$ is a prime number such that $p_i\equiv3\pmod4$ for $i\in\{1,\cdots,t+1\}$ and $p_{t+1}\nmid j$.

        \item[(v)] $B_{11}^{(11)}\left(11p_1^2\cdots p_{t+1}^2n+\frac{11\big(p_1^2\cdots p_t^2p_{t+1}(5p_{t+1}+12j)-5\big)}{12}\right)\equiv0\pmod{11},$ where $p_i\not\in\{3,11\}$ is a prime number such that $p_i\equiv3\pmod4$ for $i\in\{1,\cdots,t+1\}$ and $p_{t+1}\nmid j$.
    \end{itemize}
\end{theorem}

\begin{cor} \label{cor1}
    For any $n\geq0$ and $t\geq0$, we have
    \begin{itemize}
        \item[(i)] $B_3^{(12)}\left(3p^{2t+2}n+\left(p^{2t+1}(p+3j)-1\right)\right)\equiv0\pmod3,$ where $p$ is a prime number such that $p\equiv2\pmod3$ and $p\nmid j$.

        \item[(ii)] $B_3^{(15)}\left(3p^{2t+2}n+\frac{p^{2t+1}(5p+12j)-5}{4}\right)\equiv0\pmod3,$ where $p\neq3$ is a prime number such that $p\equiv3\pmod4$ and $p\nmid j$.

        \item[(iii)] $B_5^{(10)}\left(5p^{2t+2}n+\frac{5\big(p^{2t+1}(p+3j)-1\big)}{3}\right)\equiv0\pmod5,$ where $p\neq5$ is a prime number such that $p\equiv2\pmod3$ and $p\nmid j$.

        \item[(iv)] $B_7^{(7)}\left(7p^{2t+2}n+\frac{7\big(p^{2t+1}(p+4j)-1\big)}{4}\right)\equiv0\pmod7$, where $p\neq7$ is a prime number such that $p\equiv3\pmod4$ and $p\nmid j$.

        \item[(v)] $B_{11}^{(11)}\left(11p^{2t+2}n+\frac{11\big(p^{2t+1}(5p+12j)-5\big)}{12}\right)\equiv0\pmod{11},$ where $p\not\in\{3,11\}$ is a prime number such that $p\equiv3\pmod4$ and $p\nmid j$.
    \end{itemize}
\end{cor}

\begin{theorem} \label{thm2}
    For $n\geq0$ and $t\geq0$, we have
    \begin{itemize}
        \item[(i)] $B_5^{(6)}\left(p^{4t}n+(p^{4t}-1)\right)\equiv p^{2t}B_5^{(6)}(n)\pmod5,$ if $B_5^{(6)}\left(p-1\right)\equiv0\pmod5$, where $p\neq5$ is a prime number.

        \item[(ii)] $B_7^{(6)}\left(7p^{4t}n+\frac{7p^{4t}-3}{2}\right)\equiv p^{4t}B_7^{(6)}(7n+2)\pmod7,$ if $B_7^{(6)}\left(\frac{7p-3}{2}\right)\equiv0\pmod7$, where $p\neq7$ is an odd prime number.
    \end{itemize}
\end{theorem}

\begin{theorem} \label{thm3}
    For $n\geq0$ and $t\geq0$, we have
    \begin{itemize}
        \item[(i)] $B_3^{(12)}\left(6p^{2t+2}n+\left(p^{2t+1}(4p+3j)-1\right)\right)\equiv0\pmod6,$ where $p$ is a prime number such that $p\equiv2\pmod3$ and $j$ is an even number such that $p\nmid j$.

        \item[(ii)] $B_5^{(10)}\left(10p^{2t+2}n+\frac{5\big(p^{2t+1}(4p+3j)-1\big)}{3}\right)\equiv0\pmod{10},$ where $p\neq5$ is a prime number such that $p\equiv2\pmod3$ and $j$ is an even number such that $p\nmid j$.

        \item[(iii)] $B_3^{(15)}\left(15p^{2t+2}n+\frac{p^{2t+1}\left((12\alpha+5)p+12j\right)-5}{4}\right)\equiv0\pmod{15},$ where $\alpha\in\{1,2,3,4\}$, $p\neq3$ is a prime number such that $p\equiv3\pmod4$ and $j$ is such that $5\mid j$ but $p\nmid j$.
    \end{itemize}
\end{theorem}

\begin{theorem} \label{thm4}
    For any $n\geq0$ and $t\geq0$, we have
    \begin{itemize}
        \item[(i)] $B_5^{(20t+1)}\left(20n+14\right)\equiv0\pmod{10}$.

        \item[(ii)] $B_5^{(20t+9)}\left(20n+\alpha\right)\equiv0\pmod{10},$ where $\alpha\in\{14,18\}$.

        \item[(iii)] $B_5^{(20t+17)}\left(20n+\alpha\right)\equiv0\pmod{10},$ where $\alpha\in\{2,14,18\}$.
    \end{itemize}
\end{theorem}

\begin{theorem} \label{thm5}
   For $n\geq0$ and $t\geq0$, we have
   \begin{itemize}
       \item[(i)] $B_{35}^{(35t+4)}(35n+\alpha)\equiv0\pmod{35}$, where $\alpha\in\{13,34\}$.

       \item[(ii)] $B_{35}^{(35t+6)}(35n+\alpha)\equiv0\pmod{35}$, where $\alpha\in\{4,24,34\}$.

       \item[(iii)] $B_{35}^{(35t+11)}(35n+34)\equiv0\pmod{35}$.

       \item[(iv)] $B_{35}^{(35t+27)}(35n+\alpha)\equiv0\pmod{35}$, where $\alpha\in\{3,4,13,17,18,24,27,32,34\}$.

       \item[(v)] $B_{35}^{(35t+32)}(35n+\alpha)\equiv0\pmod{35}$, where $\alpha\in\{13,27,34\}$.

       \item[(vi)] $B_{35}^{(35t+34)}(35n+\alpha)\equiv0\pmod{35}$, where $\alpha\in\{3,4,13,18,24,34\}$.
   \end{itemize}
\end{theorem}

\begin{theorem} \label{thm6}
    For $n\geq0$ and $t\geq0$, we have
    \begin{itemize}
        \item[(i)] $B_{55}^{(55t+21)}(55n+\alpha)\equiv0\pmod{55}$, where $\alpha\in\{9,14,19,39,54\}$.

        \item[(ii)] $B_{55}^{(55t+32)}(55n+\alpha)\equiv0\pmod{55}$, where $\alpha\in\{3,8,9,14,17,19,28,32,39,42,43,47,52,53,54\}$.

        \item[(iii)] $B_{55}^{(55t+54)}(55n+\alpha)\equiv0\pmod{55}$, where $\alpha\in\{3,8,9,14,19,28,39,43,53,54\}$.
    \end{itemize}
\end{theorem}
Note that Corollary \ref{cor1} (i), (ii), (iii), (iv), (v) can be proved by taking $p_1=\cdots=p_{t+1}=p$ in Theorem \ref{thm1} (i), (ii), (iii), (iv), (v), respectively.

\section{\Large{\textbf{Preliminary}}}

In this section, we list down a few results which will be used while proving our main theorems. Before stating the first lemma, we write, for any $r\geq1$,
\begin{align} \label{eq:20}
    E_1^r=\sum_{n=0}^\infty a_r(n)q^n.
\end{align}

\begin{lemma}[\cite{newman}, Theorem 1] \label{lem1}
    Suppose that $r$ is even, $0<r\leq24$. Let $p$ be a prime such that $r(p-1)\equiv0\pmod{24}$. Set $\delta=r(p-1)/24$. Then for any $n\geq0$
    \begin{align*}
        a_r(pn+\delta)=a_r(\delta)\cdot a_r(n)-p^{\frac{r}{2}-1}\cdot a_r\left(\frac{n-\delta}{p}\right),
    \end{align*}
    with the convention that $a_r(x)=0$, if $x$ is not an integer.
\end{lemma}

Let us list down the following $2$-dissection, $5$-dissection, $7$-dissection, and $11$-dissection results.

\begin{lemma}[\cite{hirchhorn&sellers}, Theorem 2.1]
    We have,
    \begin{align} \label{eq:33}
        \frac{E_5}{E_1}=\frac{E_8E_{20}^2}{E_2^2E_{40}}+q\frac{E_4^3E_{10}E_{40}}{E_2^3E_8E_{20}}.
    \end{align}
\end{lemma}

\begin{lemma}[\cite{ramanujan} and \cite{watson1938}]
    For $R(q):=\prod_{m=1}^\infty\frac{\big(1-q^{5m-4}\big)\big(1-q^{5m-1}\big)}{\big(1-q^{5m-3}\big)\big(1-q^{5m-2}\big)}$, we have
    \begin{align} \label{eq:35}
        E_1=E_{25}\left(\frac{1}{R(q^5)}-q-q^2R\left(q^5\right)\right).
    \end{align}
\end{lemma}

\begin{lemma}[\cite{hirschhorn}, (10.5.1)]
    For $A_i(q):=\prod_{m=1}^\infty(1-q^{7m-i})(1-q^{7m-7+i})$, $i\in\{1,2,3\}$, we have
    \begin{align} \label{eq:25}
        E_1=E_{49}\left(\frac{A_2(q^7)}{A_1(q^7)}-q\frac{A_3(q^7)}{A_2(q^7)}-q^2+q^5\frac{A_1(q^7)}{A_3(q^7)}\right).
    \end{align}
\end{lemma}

\begin{lemma}[\cite{hirschhorn}, (10.6.1)]
    For $B_i(q):=\prod_{m=1}^\infty(1-q^{11m-i})(1-q^{11m-11+i})$, $i\in\{1,2,3,4,5\}$, we have
    \begin{align} \label{eq:42}
        E_1=E_{121}\left(\frac{B_4(q^{11})}{B_2(q^{11})}-q\frac{B_2(q^{11})}{B_1(q^{11})}-q^2\frac{B_5(q^{11})}{B_3(q^{11})}+q^5+q^7\frac{B_3(q^{11})}{B_4(q^{11})}-q^{15}\frac{B_1(q^{11})}{B_5(q^{11})}\right).
    \end{align}
\end{lemma}

We conclude this section by writing down Dedekind's eta-function defined by
\begin{align} \label{eq:4}
    \eta(z):=q^{1/24}\prod_{m=1}^\infty\left(1-q^m\right).
\end{align}

\section{\Large\textbf{Congruences Modulo a Prime Number}}

In this section, we will prove Theorem \ref{thm1}, and Theorem \ref{thm2}.

\begin{proof}[\textbf{Proof of Theorem \ref{thm1}}]
\textbf{{Part (i)}:} Let us take $\ell=3$ and $r=12$ in (\ref{eq:2}) to get
\begin{align*}
    \sum_{n=0}^\infty B_3^{(12)}(n)q^n=\frac{E_3^{12}}{E_1^{12}}.
\end{align*}
Now, by (\ref{eq:3}), we get
\begin{align*}
    \sum_{n=0}^\infty B_3^{(12)}(n)q^n\equiv\frac{E_3^{12}}{E_3^{4}}\equiv E_3^8\pmod3.
\end{align*}
In view of (\ref{eq:4}), we get
\begin{align*}
    \sum_{n=0}^\infty B_3^{(12)}(3n)q^{3n+1}\equiv\eta^8(3z)\pmod3.
\end{align*}
Writing $\eta^8(3z)=\sum_{n=0}^\infty a(n)q^n$, we get
\begin{align} \label{eq:6}
    B_3^{(12)}(3n)\equiv a(3n+1)\pmod3.
\end{align}
Moreover, one can verify that $a(n)=0$, if $n\not\equiv1\pmod3$.

Serre \cite{Serre} showed that $\eta^8(3z)$ is a normalized Hecke eigenform of weight $4$, level $9$. Therefore, for any prime number $p$, operating the Hecke operator $T_p$ on $\eta^8(3z)$, we get
\begin{align} \label{eq:5}
    \eta^8(3z)\mid T_p=a(p)\cdot\eta^8(3z).
\end{align}
Since
\begin{align*}
   \eta^8(3z)\mid T_p=&\sum_{n=1}^\infty\left(a(pn)+p^{4-1}\cdot a\left(\frac{n}{p}\right)\right)q^n\\
   =&\sum_{n=1}^\infty\left(a(pn)+p^3\cdot a\left(\frac{n}{p}\right)\right)q^n,
\end{align*}
therefore, by comparing coefficients on both sides of (\ref{eq:5}), we get
\begin{align*}
    a(pn)+p^3\cdot a\left(\frac{n}{p}\right)=a(p)\cdot a(n)\quad\forall\quad n\geq1.
\end{align*}

Let $p$ be a prime number such that $p\equiv2\pmod3$; then we have
\begin{align} \label{eq:11}
    a(pn)+p^3\cdot a\left(\frac{n}{p}\right)=0.
\end{align}
Substituting $n$ by $pn+r$ with $p\nmid r$ in the above equation, we get $$a(p^2n+pr)=0.$$ Now, substituting $n$ by $3n+1-pr$, we get $$a(3p^2n+p^2-p^3r+pr)=0.$$ But, $3p^2n+p^2-p^3r+pr=3\left(p^2n+\frac{p^2-1}{3}+\frac{pr(1-p^2)}{3}\right)+1$. Thus, by (\ref{eq:6}), we get
\begin{align*}
    B_3^{(12)}\left(3\left(p^2n+\frac{p^2-1}{3}+\frac{pr(1-p^2)}{3}\right)\right)\equiv0\pmod3.
\end{align*}
As $gcd\left(\frac{1-p^2}{3},p\right)=1$, so we have $p\nmid r$ if and only if $p\nmid j$, where $j:=\frac{r(1-p^2)}{3}$. Hence, for $j\not\equiv0\pmod p$, we get
\begin{align} \label{eq:8}
    B_3^{(12)}\left(3p^2n+p^2-1+3pj\right)\equiv0\pmod3.
\end{align}

Substituting $n$ by $pn$ in (\ref{eq:11}), we get
\begin{align*}
    a(p^2n)=-p^3\cdot a(n)\equiv a(n)\pmod3.
\end{align*}
Now, substituting $n$ by $3n+1$, we get
\begin{align*}
    a(3p^2n+p^2)\equiv a(3n+1)\pmod3.
\end{align*}
But, $3p^2n+p^2=3\left(p^2n+\frac{p^2-1}{3}\right)+1$. Thus, by (\ref{eq:6}), we get
\begin{align} \label{eq:7}
    B_3^{(12)}\left(3\left(p^2n+\frac{p^2-1}{3}\right)\right)\equiv B_3^{(12)}(3n)\pmod3.
\end{align}
Let $t\geq0$ and $p_i$ be a prime number such that $p_i\equiv2\pmod3$ for $i\in\{1,\cdots,t+1\}$. Since
\begin{align*}
    p_1^2\cdots p_t^2n+\frac{p_1^2\cdots p_t^2-1}{3}=p_1^2\left(p_2^2\cdots p_t^2n+\frac{p_2^2\cdots p_t^2-1}{3}\right)+\frac{p_1^2-1}{3},
\end{align*}
therefore by repeated use of (\ref{eq:7}), we get
\begin{align*}
    B_3^{(12)}\left(3p_1^2\cdots p_t^2n+(p_1^2\cdots p_t^2-1)\right)\equiv B_3^{(12)}(3n)\pmod3.
\end{align*}
Substituting $n$ by $p_{t+1}^2n+\frac{p_{t+1}^2-1}{3}+p_{t+1}j$ with $p_{t+1}\nmid j$, we get the desired result in view of (\ref{eq:8}).

\textbf{{Part (ii)}:} Let us take $\ell=3$ and $r=15$ in (\ref{eq:2}) to get
\begin{align*}
    \sum_{n=0}^\infty B_3^{(15)}(n)q^n=\frac{E_3^{15}}{E_1^{15}}.
\end{align*}
Now, by (\ref{eq:3}), we get
\begin{align*}
    \sum_{n=0}^\infty B_3^{(15)}(n)q^n\equiv\frac{E_3^{15}}{E_3^{5}}\equiv E_3^{10}\pmod3.
\end{align*}
In view of (\ref{eq:4}), we get
\begin{align*}
    \sum_{n=0}^\infty B_3^{(15)}(3n)q^{12n+5}\equiv\eta^{10}(12z)\pmod3.
\end{align*}
Writing $\eta^{10}(12z)=\sum_{n=0}^\infty a(n)q^n$, we get
\begin{align} \label{eq:10}
    B_3^{(15)}(3n)\equiv a(12n+5)\pmod3.
\end{align}

Serre \cite{Serre} showed that $\eta^{10}(12z)=\frac{1}{96}(\varphi_+(z)-\varphi_-(z))$, where $\varphi_+(z)$ and $\varphi_-(z)$ are normalized Hecke eigenforms of weight $5$, level $2^43^2$ and character $\chi$ defined as
\begin{align*}
    \chi(n)=\begin{cases}
    (-1)^{\frac{n-1}{2}}&\text{ if $gcd(n,2^43^2)=1$},\\
    0&\text{ otherwise}.
    \end{cases}
\end{align*}
If we write $\varphi_\pm(z)=\sum_{n=0}^\infty b_\pm(n)q^n$, then for any prime number $p$, operating the Hecke operator $T_p$ on $\eta^{10}(12z)$, we get
\begin{align} \label{eq:9}
    \eta^{10}(12z)\mid T_p=\frac{1}{96}\left(\varphi_+(z)\mid T_p-\varphi_-(z)\mid T_p\right)=\frac{1}{96}\left(b_+(p)\cdot\varphi_+(z)-b_-(p)\cdot\varphi_-(z)\right).
\end{align}
Since
\begin{align*}
   \eta^{10}(12z)\mid T_p=&\sum_{n=1}^\infty\left(a(pn)+\chi(p)\cdot p^{5-1}\cdot a\left(\frac{n}{p}\right)\right)q^n\\
   =&\sum_{n=1}^\infty\left(a(pn)+(-1)^{\frac{p-1}{2}}p^4\cdot a\left(\frac{n}{p}\right)\right)q^n,
\end{align*}
therefore, by comparing coefficients on both sides of (\ref{eq:9}), we get
\begin{align*}
    a(pn)+(-1)^{\frac{p-1}{2}}p^4\cdot a\left(\frac{n}{p}\right)=\frac{1}{96}\left(b_+(p)\cdot b_+(n)-b_-(p)\cdot b_-(n)\right)\quad\forall\quad n\geq1.
\end{align*}
Moreover, Serre \cite{Serre} further established that $b_\pm(n)=0$, if $n\not\equiv1\pmod4$.

Let $p\neq3$ be a prime number such that $p\equiv3\pmod4$; then we have
\begin{align} \label{eq:12}
    a(pn)-p^4\cdot a\left(\frac{n}{p}\right)=0.
\end{align}
Substituting $n$ by $pn+r$ with $p\nmid r$ in the above equation, we get $$a(p^2n+pr)=0.$$ Now, substituting $n$ by $12n+5-pr$, we get $$a(12p^2n+5p^2-p^3r+pr)=0.$$ But, $12p^2n+5p^2-p^3r+pr=12\left(p^2n+\frac{5(p^2-1)}{12}+\frac{pr(1-p^2)}{12}\right)+5$. Thus, by (\ref{eq:10}), we get
\begin{align*}
    B_{3}^{(15)}\left(3\left(p^2n+\frac{5(p^2-1)}{12}+\frac{pr(1-p^2)}{12}\right)\right)\equiv0\pmod{3}.
\end{align*}
As $gcd\left(\frac{1-p^2}{12},p\right)=1$, so we have $p\nmid r$ if and only if $p\nmid j$, where $j:=\frac{r(1-p^2)}{12}$. Hence, for $j\not\equiv0\pmod p$, we get
\begin{align} \label{eq:14}
    B_{3}^{(15)}\left(3p^2n+\frac{5(p^2-1)}{4}+3pj\right)\equiv0\pmod{3}.
\end{align}

Substituting $n$ by $pn$ in (\ref{eq:12}), we get
\begin{align*}
    a(p^2n)=p^4\cdot a(n)\equiv a(n)\pmod3.
\end{align*}
Now, substituting $n$ by $12n+5$, we get
\begin{align*}
    a(12p^2n+5p^2)\equiv a(12n+5)\pmod3.
\end{align*}
But, $12p^2n+5p^2=12\left(p^2n+\frac{5(p^2-1)}{12}\right)+5$. Thus, by (\ref{eq:10}), we get
\begin{align} \label{eq:13}
    B_{3}^{(15)}\left(3\left(p^2n+\frac{5(p^2-1)}{12}\right)\right)\equiv B_{3}^{(15)}(3n)\pmod{3}.
\end{align}
Let $t\geq0$ and $p_i\neq3$ be a prime number such that $p_i\equiv3\pmod4$ for $i\in\{1,\cdots,t+1\}$. Since
\begin{align*}
    p_1^2\cdots p_t^2n+\frac{5(p_1^2\cdots p_t^2-1)}{12}=p_1^2\left(p_2^2\cdots p_t^2n+\frac{5(p_2^2\cdots p_t^2-1)}{12}\right)+\frac{5(p_1^2-1)}{12},
\end{align*}
therefore by repeated use of (\ref{eq:13}), we get
\begin{align*}
    B_{3}^{(15)}\left(3p_1^2\cdots p_t^2n+\frac{5(p_1^2\cdots p_t^2-1)}{4}\right)\equiv B_{3}^{(15)}(3n)\pmod{3}.
\end{align*}
Substituting $n$ by $p_{t+1}^2n+\frac{5(p_{t+1}^2-1)}{12}+p_{t+1}j$ with $p_{t+1}\nmid j$, we get the desired result in view of (\ref{eq:14}).

\textbf{{Part (iii)}:} Let us take $\ell=5$ and $r=10$ in (\ref{eq:2}) to get
\begin{align*}
    \sum_{n=0}^\infty B_5^{(10)}(n)q^n=\frac{E_5^{10}}{E_1^{10}}.
\end{align*}
Now, by (\ref{eq:3}), we get
\begin{align*}
    \sum_{n=0}^\infty B_5^{(10)}(n)q^n\equiv\frac{E_5^{10}}{E_5^{2}}\equiv E_5^8\pmod5.
\end{align*}
In view of (\ref{eq:4}), we get
\begin{align*}
    \sum_{n=0}^\infty B_5^{(10)}(5n)q^{3n+1}\equiv\eta^8(3z)\pmod5.
\end{align*}
Rest of the proof follows similarly as \textbf{Part (i)}.

\textbf{{Part (iv)}:} Let us take $\ell=7$ and $r=7$ in (\ref{eq:2}) to get
\begin{align*}
    \sum_{n=0}^\infty B_7^{(7)}(n)q^n=\frac{E_7^7}{E_1^7}.
\end{align*}
Now, by (\ref{eq:3}), we get
\begin{align*}
    \sum_{n=0}^\infty B_7^{(7)}(n)q^n\equiv\frac{E_7^7}{E_7}\equiv E_7^6\pmod7.
\end{align*}
In view of (\ref{eq:4}), we get
\begin{align*}
    \sum_{n=0}^\infty B_7^{(7)}(7n)q^{4n+1}\equiv\eta^6(4z)\pmod7.
\end{align*}
Writing $\eta^6(4z)=\sum_{n=0}^\infty a(n)q^n$, we get
\begin{align} \label{eq:16}
    B_7^{(7)}(7n)\equiv a(4n+1)\pmod7.
\end{align}
Moreover, one can verify that $a(n)=0$, if $n\not\equiv1\pmod4$.

Serre \cite{Serre} showed that $\eta^6(4z)$ is a normalized Hecke eigenform of weight $3$, level $16$ and character $\chi$ defined as
\begin{align*}
    \chi(n)=\begin{cases}
    (-1)^{\frac{n-1}{2}}&\text{ if $n$ is odd},\\
    0&\text{ if $n$ is even}.
    \end{cases}
\end{align*}
Therefore, for any prime number $p$, operating the Hecke operator $T_p$ on $\eta^6(4z)$, we get
\begin{align} \label{eq:15}
    \eta^6(4z)\mid T_p=a(p)\cdot\eta^6(4z).
\end{align}
Since
\begin{align*}
   \eta^6(4z)\mid T_p=&\sum_{n=1}^\infty\left(a(pn)+\chi(p)\cdot p^{3-1}\cdot a\left(\frac{n}{p}\right)\right)q^n\\
   =&\sum_{n=1}^\infty\left(a(pn)+(-1)^{\frac{p-1}{2}}p^2\cdot a\left(\frac{n}{p}\right)\right)q^n,
\end{align*}
therefore, by comparing coefficients on both sides of (\ref{eq:15}), we get
\begin{align*}
    a(pn)+(-1)^{\frac{p-1}{2}}p^2\cdot a\left(\frac{n}{p}\right)=a(p)\cdot a(n)\quad\forall\quad n\geq1.
\end{align*}

Let $p\neq7$ be a prime number such that $p\equiv3\pmod4$; then we have
\begin{align} \label{eq:17}
    a(pn)-p^2\cdot a\left(\frac{n}{p}\right)=0.
\end{align}
Substituting $n$ by $pn+r$ with $p\nmid r$ in the above equation, we get $$a(p^2n+pr)=0.$$ Now, substituting $n$ by $4n+1-pr$, we get $$a(4p^2n+p^2-p^3r+pr)=0.$$ But, $4p^2n+p^2-p^3r+pr=4\left(p^2n+\frac{p^2-1}{4}+\frac{pr(1-p^2)}{4}\right)+1$. Thus, by (\ref{eq:16}), we get
\begin{align*}
    B_7^{(7)}\left(7\left(p^2n+\frac{p^2-1}{4}+\frac{pr(1-p^2)}{4}\right)\right)\equiv0\pmod7.
\end{align*}
As $gcd\left(\frac{1-p^2}{4},p\right)=1$, so we have $p\nmid r$ if and only if $p\nmid j$, where $j:=\frac{r(1-p^2)}{4}$. Hence, for $j\not\equiv0\pmod p$, we get
\begin{align} \label{eq:19}
    B_7^{(7)}\left(7p^2n+\frac{7(p^2-1)}{4}+7pj\right)\equiv0\pmod7.
\end{align}

Substituting $n$ by $pn$ in (\ref{eq:17}), we get
\begin{align*}
    a(p^2n)=p^2\cdot a(n).
\end{align*}
Now, substituting $n$ by $4n+1$, we get
\begin{align*}
    a(4p^2n+p^2)=p^2\cdot a(4n+1).
\end{align*}
But, $4p^2n+p^2=4\left(p^2n+\frac{p^2-1}{4}\right)+1$. Thus, by (\ref{eq:16}), we get
\begin{align} \label{eq:18}
    B_7^{(7)}\left(7\left(p^2n+\frac{p^2-1}{4}\right)\right)\equiv p^2\cdot B_7^{(7)}(7n)\pmod7.
\end{align}
Let $t\geq0$ and $p_i\neq7$ be a prime number such that $p_i\equiv3\pmod4$ for $i\in\{1,\cdots,t+1\}$. Since
\begin{align*}
    p_1^2\cdots p_t^2n+\frac{p_1^2\cdots p_t^2-1}{4}=p_1^2\left(p_2^2\cdots p_t^2n+\frac{p_2^2\cdots p_t^2-1}{4}\right)+\frac{p_1^2-1}{4},
\end{align*}
therefore by repeated use of (\ref{eq:18}), we get
\begin{align*}
    B_7^{(7)}\left(7p_1^2\cdots p_t^2n+\frac{7(p_1^2\cdots p_t^2-1)}{4}\right)\equiv(p_1^2\cdots p_{t+1}^2)\cdot B_7^{(7)}(7n)\pmod7.
\end{align*}
Substituting $n$ by $p_{t+1}^2n+\frac{p_{t+1}^2-1}{4}+p_{t+1}j$ with $p_{t+1}\nmid j$, we get the desired result in view of (\ref{eq:19}).

\textbf{{Part (v)}:} Let us take $\ell=11$ and $r=11$ in (\ref{eq:2}) to get
\begin{align*}
    \sum_{n=0}^\infty B_{11}^{(11)}(n)q^n=\frac{E_{11}^{11}}{E_1^{11}}.
\end{align*}
Now, by (\ref{eq:3}), we get
\begin{align*}
    \sum_{n=0}^\infty B_{11}^{(11)}(n)q^n\equiv\frac{E_{11}^{11}}{E_{11}}\equiv E_{11}^{10}\pmod{11}.
\end{align*}
In view of (\ref{eq:4}), we get
\begin{align*}
    \sum_{n=0}^\infty B_{11}^{(11)}(11n)q^{12n+5}\equiv\eta^{10}(12z)\pmod{11}.
\end{align*}
Rest of the proof follows similarly as \textbf{Part (ii)}.
\end{proof}

\begin{proof}[\textbf{Proof of Theorem \ref{thm2}}]
\textbf{{Part (i)}:} Let us take $\ell=5$ and $r=6$ in (\ref{eq:2}) to get
\begin{align*}
    \sum_{n=0}^\infty B_5^{(6)}(n)q^n=\frac{E_5^{6}}{E_1^{6}}.
\end{align*}
Now, by (\ref{eq:3}), we get
\begin{align*}
    \sum_{n=0}^\infty B_5^{(6)}(n)q^n\equiv\frac{E_1^{30}}{E_1^{6}}\equiv E_1^{24}\pmod{5}.
\end{align*}
In view of (\ref{eq:20}), we get
\begin{align} \label{eq:24}
    B_5^{(6)}(n)\equiv a_{24}(n)\pmod5 \quad \forall \quad n\geq0.
\end{align}

Let us apply Lemma \ref{lem1} taking $r=24$ to get that for a prime $p$,
\begin{align*}
    a_{24}\left(pn+(p-1)\right)=a_{24}\left(p-1\right)\cdot a_{24}(n)-p^{11}\cdot a_{24}\left(\frac{n-(p-1)}{p}\right).
\end{align*}
Substituting $n$ by $pn+(p-1)$ repeatedly, we get
\begin{align}
    a_{24}\left(p^2n+(p^2-1)\right)&=a_{24}\left(p-1\right)\cdot a_{24}\left(pn+(p-1)\right)-p^{11}\cdot a_{24}\left(n\right), \label{eq:21}\\
    a_{24}\left(p^3n+(p^3-1)\right)&=a_{24}\left(p-1\right)\cdot a_{24}\left(p^2n+(p^2-1)\right)-p^{11}\cdot a_{24}\left(pn+(p-1)\right), \label{eq:22}\\
    a_{24}\left(p^4n+(p^4-1)\right)&=a_{24}\left(p-1\right)\cdot a_{24}\left(p^3n+(p^3-1)\right)-p^{11}\cdot a_{24}\left(p^2n+(p^2-1)\right). \label{eq:23}
\end{align}
In (\ref{eq:23}), replacing $a_{24}\left(p^3n+(p^3-1)\right)$ and $a_{24}\left(p^2n+(p^2-1)\right)$ according to (\ref{eq:21}) and (\ref{eq:22}), respectively, we get after simplifying
\begin{align*}
    a_{24}\left(p^4n+(p^4-1)\right)&=a_{24}\left(p-1\right)\cdot\left(a_{24}\left(p-1\right)^2-2p^{11}\right)\cdot a_{24}\left(pn+(p-1)\right)\\
    &\quad-p^{11}\cdot\left(a_{24}\left(p-1\right)^2-p^{11}\right)\cdot a_{24}(n).
\end{align*}
Now, by (\ref{eq:24}), we get
\begin{align*}
    B_{5}^{(6)}\left(p^4n+(p^4-1)\right)&\equiv B_{5}^{(6)}\left(p-1\right)\cdot\left(B_{5}^{(6)}\left(p-1\right)^2-2p^{11}\right)\cdot B_{5}^{(6)}\left(pn+(p-1)\right)\\
    &\quad-p^{11}\cdot\left(B_{5}^{(6)}\left(p-1\right)^2-p^{11}\right)\cdot B_{5}^{(6)}(n)\pmod5.
\end{align*}
Let us assume that $p\neq5$ and that $B_5^{(6)}\left(p-1\right)\equiv0\pmod5$ to get
\begin{align*}
    B_5^{(6)}\left(p^4n+(p^4-1)\right)\equiv p^{2}B_5^{(6)}(n)\pmod5.
\end{align*}
We get the desired result by inductively substituting $n$ by $p^4n+(p^4-1)$ in the above congruence.

\textbf{{Part (ii)}:} Let us take $\ell=7$ and $r=6$ in (\ref{eq:2}) to get
\begin{align*}
    \sum_{n=0}^\infty B_7^{(6)}(n)q^n=\frac{E_7^{6}}{E_1^{6}}.
\end{align*}
The above equation can be rewritten as
\begin{align*}
    \sum_{n=0}^\infty B_7^{(6)}(n)q^n=\frac{E_7^{6}}{E_1^{7}}\cdot E_1.
\end{align*}
Now, by (\ref{eq:25}), we get
\begin{align*}
    \sum_{n=0}^\infty B_7^{(6)}(n)q^n=\frac{E_7^{6}}{E_1^{7}}\cdot E_{49}\Bigg(\frac{B\big(q^7\big)}{C(q^7)}-q\frac{A\big(q^7\big)}{B(q^7)}-q^2+q^5\frac{C\big(q^7\big)}{A(q^7)}\Bigg).
\end{align*}
In view of (\ref{eq:3}), we get
\begin{align*}
    \sum_{n=0}^\infty B_7^{(6)}(n)q^n\equiv E_7^{12}\Bigg(\frac{B\big(q^7\big)}{C(q^7)}-q\frac{A\big(q^7\big)}{B(q^7)}-q^2+q^5\frac{C\big(q^7\big)}{A(q^7)}\Bigg)\pmod7.
\end{align*}

Collecting the terms involving $q^{7n+2}$, we get
\begin{align*}
    \sum_{n=0}^\infty B_7^{(6)}(7n+2)q^{7n+2}\equiv6q^2E_7^{12}\pmod7.
\end{align*}
Now, dividing both sides by $q^2$ and substituting $q$ by $q^{1/7}$, we get
\begin{align*}
    \sum_{n=0}^\infty B_7^{(6)}(7n+2)q^n\equiv6E_1^{12}\pmod7.
\end{align*}
In view of (\ref{eq:20}), we get
\begin{align} \label{eq:29}
    B_7^{(6)}(7n+2)\equiv 6a_{12}(n)\pmod7 \quad \forall \quad n\geq0.
\end{align}

Let us apply Lemma \ref{lem1} taking $r=12$ to get that for an odd prime $p$,
\begin{align*}
    a_{12}\left(pn+\frac{p-1}{2}\right)=a_{12}\left(\frac{p-1}{2}\right)\cdot a_{12}(n)-p^{5}\cdot a_{12}\left(\frac{n-\frac{p-1}{2}}{p}\right).
\end{align*}
Substituting $n$ by $pn+\frac{p-1}{2}$ repeatedly, we get
\begin{align}
    a_{12}\left(p^2n+\frac{p^2-1}{2}\right)&=a_{12}\left(\frac{p-1}{2}\right)\cdot a_{12}\left(pn+\frac{p-1}{2}\right)-p^{5}\cdot a_{12}\left(n\right), \label{eq:26}\\
    a_{12}\left(p^3n+\frac{p^3-1}{2}\right)&=a_{12}\left(\frac{p-1}{2}\right)\cdot a_{12}\left(p^2n+\frac{p^2-1}{2}\right)-p^{5}\cdot a_{12}\left(pn+\frac{p-1}{2}\right), \label{eq:27}\\
    a_{12}\left(p^4n+\frac{p^4-1}{2}\right)&=a_{12}\left(\frac{p-1}{2}\right)\cdot a_{12}\left(p^3n+\frac{p^3-1}{2}\right)-p^{5}\cdot a_{12}\left(p^2n+\frac{p^2-1}{2}\right). \label{eq:28}
\end{align}
In (\ref{eq:28}), replacing $a_{12}\left(p^3n+\frac{p^3-1}{2}\right)$ and $a_{12}\left(p^2n+\frac{p^2-1}{2}\right)$ according to (\ref{eq:26}) and (\ref{eq:27}), respectively, we get after simplifying
\begin{align*}
    a_{12}\left(p^4n+\frac{p^4-1}{2}\right)&=a_{12}\left(\frac{p-1}{2}\right)\cdot\left(a_{12}\left(\frac{p-1}{2}\right)^2-2p^5\right)\cdot a_{12}\left(pn+\frac{p-1}{2}\right)\\
    &\quad-p^5\cdot\left(a_{12}\left(\frac{p-1}{2}\right)^2-p^5\right)\cdot a_{12}(n).
\end{align*}
Now, by (\ref{eq:29}), we get
\begin{align*}
    6B_7^{(6)}\left(7p^4n+\frac{7p^4-3}{2}\right)&\equiv6B_7^{(6)}\left(\frac{7p-3}{2}\right)\cdot\left(6B_7^{(6)}\left(\frac{7p-3}{2}\right)^2-2p^5\right)\cdot 6B_7^{(6)}\left(7pn+\frac{7p-3}{2}\right)\\
    &\quad-p^5\cdot\left(6B_7^{(6)}\left(\frac{7p-3}{2}\right)^2-p^5\right)\cdot 6B_7^{(6)}(7n+2)\pmod7.
\end{align*}
Let us assume that $p\neq7$ and that $B_7^{(6)}\left(\frac{7p-3)}{2}\right)\equiv0\pmod7$. Then we get
\begin{align*}
    B_7^{(6)}\left(7p^4n+\frac{7p^4-3}{2}\right)\equiv p^{4}B_7^{(6)}(7n+2)\pmod7.
\end{align*}
We get the desired result by inductively substituting $n$ by $p^4n+\frac{p^4-1}{2}$ in the above congruence.
\end{proof}

\section{\textbf{Congruences Modulo a Composite Number}}

In this section, we will prove Theorem \ref{thm3}, Theorem \ref{thm4}, Theorem \ref{thm5}, and Theorem \ref{thm6}.

\begin{proof} [\textbf{Proof of Theorem \ref{thm3}}]
\textbf{{Part (i)}:} Let us take $\ell=3$ and $r=12$ in (\ref{eq:2}) to get
\begin{align*}
    \sum_{n=0}^\infty B_3^{(12)}(n)q^n=\frac{E_3^{12}}{E_1^{12}}.
\end{align*}
Now, by (\ref{eq:3}), we get
\begin{align*}
    \sum_{n=0}^\infty B_3^{(12)}(n)q^n\equiv\frac{E_6^{6}}{E_2^{6}}\pmod2.
\end{align*}
Comparing coefficients of $q^{2n+1}$ on both sides, we get
\begin{align} \label{eq:30}
    B_3^{(12)}(2n+1)\equiv0\pmod2 \quad \forall \quad n\geq0.
\end{align}

On the other hand, let $p$ be a prime number such that $p\equiv2\pmod3$ and $j$ be an even number such that $p\nmid j$. By Corollary \ref{cor1} (i), we get that for any $n\geq0$ and $t\geq0$,
\begin{align*}
    B_3^{(12)}\left(3p^{2t+2}n+\left(p^{2t+1}(p+3j)-1\right)\right)\equiv0\pmod3.
\end{align*}
Substituting $n$ by $2n+1$, we get
\begin{align*}
    B_3^{(12)}\left(6p^{2t+2}n+\left(p^{2t+1}(4p+3j)-1\right)\right)\equiv0\pmod3.
\end{align*}
Note that $6p^{2t+2}n+p^{2t+1}(4p+3j)-1=2\left(3p^{2t+2}n+\frac{p^{2t+1}(4p+3j)-2}{2}\right)+1$. Thus, in view of (\ref{eq:30}), we get
\begin{align*}
    B_3^{(12)}\left(6p^{2t+2}n+\left(p^{2t+1}(4p+3j)-1\right)\right)\equiv0\pmod2,
\end{align*}
and hence, the desired result.

\textbf{{Part (ii)}:} Let us take $\ell=5$ and $r=10$ in (\ref{eq:2}) to get
\begin{align*}
    \sum_{n=0}^\infty B_5^{(10)}(n)q^n=\frac{E_5^{10}}{E_1^{10}}.
\end{align*}
Now, by (\ref{eq:3}), we get
\begin{align*}
    \sum_{n=0}^\infty B_5^{(10)}(n)q^n\equiv\frac{E_{10}^{5}}{E_2^{5}}\pmod{2}.
\end{align*}
Comparing coefficients of $q^{2n+1}$ on both sides, we get
\begin{align} \label{eq:31}
    B_5^{(10)}(2n+1)\equiv0\pmod2 \quad \forall \quad n\geq0.
\end{align}

On the other hand, let $p\neq5$ be a prime number such that $p\equiv2\pmod3$ and $j$ be an even number such that $p\nmid j$. By Corollary \ref{cor1} (iii), we get that for any $n\geq0$ and $t\geq0$,
\begin{align*}
    B_5^{(10)}\left(5p^{2t+2}n+\frac{5\big(p^{2t+1}(p+3j)-1\big)}{3}\right)\equiv0\pmod5.
\end{align*}
Substituting $n$ by $2n+1$, we get
\begin{align*}
    B_5^{(10)}\left(10p^{2t+2}n+\frac{5\big(p^{2t+1}(4p+3j)-1\big)}{3}\right)\equiv0\pmod5.
\end{align*}
Note that $10p^{2t+2}n+\frac{5\big(p^{2t+1}(4p+3j)-1\big)}{3}=2\left(5p^{2t+2}n+\frac{5p^{2t+1}(4p+3j)-8}{6}\right)+1$. Thus, in view of (\ref{eq:31}), we get
\begin{align*}
    B_5^{(10)}\left(10p^{2t+2}n+\frac{5\big(p^{2t+1}(4p+3j)-1\big)}{3}\right)\equiv0\pmod2,
\end{align*}
and hence, the desired result.

\textbf{{Part (iii)}:} Let us take $\ell=3$ and $r=15$ in (\ref{eq:2}) to get
\begin{align*}
    \sum_{n=0}^\infty B_3^{(15)}(n)q^n=\frac{E_3^{15}}{E_1^{15}}.
\end{align*}
Now, by (\ref{eq:3}), we get
\begin{align*}
    \sum_{n=0}^\infty B_3^{(15)}(n)q^n\equiv\frac{E_{15}^{3}}{E_5^{3}}\pmod{5}.
\end{align*}
Comparing coefficients of $q^{5n+\alpha}$ on both sides, where $\alpha\in\{1,2,3,4\}$, we get
\begin{align} \label{eq:32}
    B_3^{(15)}(5n+\alpha)\equiv0\pmod5 \quad \forall \quad n\geq0.
\end{align}

On the other hand, let $p\neq3$ be a prime number such that $p\equiv2\pmod3$ and $j$ be a positive integer such that $5\mid j$ but $p\nmid j$. By Corollary \ref{cor1} (ii), we get that for $n\geq0$ and $t\geq0$,
\begin{align*}
    B_3^{(15)}\left(3p^{2t+2}n+\frac{p^{2t+1}(5p+12j)-5}{4}\right)\equiv0\pmod3.
\end{align*}
Substituting $n$ by $5n+\alpha$, we get
\begin{align*}
    B_3^{(15)}\left(15p^{2t+2}n+\frac{p^{2t+1}\left((12\alpha+5)p+12j\right)-5}{4}\right)\equiv0\pmod3.
\end{align*}
Note that $15p^{2t+2}n+\frac{p^{2t+1}\left((12\alpha+5)p+12j\right)-5}{4}=5\left(3p^{2t+2}n+\frac{3p^{2t+2}\alpha-\alpha'+3p^{2t+1}j}{5}+\frac{p^{2t+2}-1}{4}\right)+\alpha'$, where $\alpha'\in\{1,2,3,4\}$ is such that $3p^{2t+2}\alpha\equiv\alpha'\pmod5$. Thus, in view of (\ref{eq:32}), we get
\begin{align*}
    B_3^{(15)}\left(15p^{2t+2}n+\frac{p^{2t+1}\left((12\alpha+5)p+12j\right)-5}{4}\right)\equiv0\pmod5,
\end{align*}
and hence, the desired result.
\end{proof}

\begin{proof}[\textbf{Proof of Theorem \ref{thm4}}]
\textbf{{Part (i)}:} Let us take $\ell=5$ and $r=20t+1$, for $t\geq0$, in (\ref{eq:2}) to get
\begin{align} \label{eq:34}
    \sum_{n=0}^\infty B_5^{(20t+1)}(n)q^n=\frac{E_5^{20t+1}}{E_1^{20t+1}}.
\end{align}
The above equation can be rewritten as
\begin{align*}
    \sum_{n=0}^\infty B_5^{(20t+1)}(n)q^n=\frac{E_5^{20t}}{E_1^{20t}}\cdot\frac{E_5}{E_1}.
\end{align*}
Now, by (\ref{eq:33}), we get
\begin{align*}
    \sum_{n=0}^\infty B_5^{(20t+1)}(n)q^n=\frac{E_5^{20t}}{E_1^{20t}}\left(\frac{E_8E_{20}^2}{E_2^2E_{40}}+q\frac{E_4^3E_{10}E_{40}}{E_2^3E_8E_{20}}\right).
\end{align*}
In view of (\ref{eq:3}), we get
\begin{align*}
    \sum_{n=0}^\infty B_5^{(20t+1)}(n)q^n\equiv\frac{E_{20}^{5t}}{E_4^{5t}}\left(\frac{E_8E_{20}^2}{E_4E_{40}}+q\frac{E_4^3E_{10}E_{40}}{E_2^3E_8E_{20}}\right)\pmod2.
\end{align*}

Collecting the terms involving $q^{2n}$, we get
\begin{align*}
    \sum_{n=0}^\infty B_5^{(20t+1)}(2n)q^{2n}\equiv\frac{E_{20}^{5t}}{E_4^{5t}}\cdot\frac{E_8E_{20}^2}{E_4E_{40}}\pmod2.
\end{align*}
Now, substituting $q$ by $q^{1/2}$, we get
\begin{align*}
    \sum_{n=0}^\infty B_5^{(20t+1)}(2n)q^{n}\equiv\frac{E_{10}^{5t}}{E_2^{5t}}\cdot\frac{E_4E_{10}^2}{E_2E_{20}}\pmod2.
\end{align*}
Comparing coefficients of $q^{2n+1}$ on both sides, we get
\begin{align} \label{eq:36}
    B_5^{(20t+1)}(4n+2)\equiv0\pmod2 \quad \forall \quad n\geq0.
\end{align}

On the other hand, (\ref{eq:34}) can be also rewritten as
\begin{align*}
    \sum_{n=0}^\infty B_5^{(20t+1)}(n)q^n=\frac{E_5^{20t+1}}{E_1^{20t+5}}\cdot E_1^4.
\end{align*}
Now, by (\ref{eq:35}), we get
\begin{align*}
    \sum_{n=0}^\infty B_5^{(20t+1)}(n)q^n=\frac{E_5^{20t+1}}{E_1^{20t+5}}\cdot E_{25}^4\left(\frac{1}{R(q^5)}-q-q^2R\left(q^5\right)\right)^4.
\end{align*}
In view of (\ref{eq:3}), we get
\begin{align*}
    \sum_{n=0}^\infty B_5^{(20t+1)}(n)q^n\equiv\frac{E_{5}^{20t+1}E_{25}^4}{E_5^{4t+1}}\Bigg(&\frac{1}{R\left(q^5\right)^4}+\frac{q}{R\left(q^5\right)^3}+\frac{2q^2}{R\left(q^5\right)^2}+\frac{3q^3}{R\left(q^5\right)}+2q^5R\left(q^5\right)\\
    &+2q^6R\left(q^5\right)^2+4q^7R\left(q^5\right)^3+q^8R\left(q^5\right)^4\Bigg)\pmod5.
\end{align*}
Comparing coefficients of $q^{5n+4}$ on both sides, we get
\begin{align} \label{eq:37}
    B_5^{(20t+1)}(5n+4)\equiv0\pmod5 \quad \forall \quad n\geq0.
\end{align}
Combining (\ref{eq:36}) and (\ref{eq:37}), we get the desired result.

\textbf{{Part (ii)}:} Let us take $\ell=5$ and $r=20t+9$, for $t\geq0$, in (\ref{eq:2}) to get
\begin{align*}
    \sum_{n=0}^\infty B_5^{(20t+9)}(n)q^n=\frac{E_5^{20t+9}}{E_1^{20t+9}}.
\end{align*}
The above equation can be rewritten as
\begin{align*}
    \sum_{n=0}^\infty B_5^{(20t+9)}(n)q^n=\frac{E_5^{20t+8}}{E_1^{20t+8}}\cdot\frac{E_5}{E_1}
\end{align*}
as well as
\begin{align*}
    \sum_{n=0}^\infty B_5^{(20t+9)}(n)q^n=\frac{E_5^{20t+9}}{E_1^{20t+10}}\cdot E_1.
\end{align*}
Following similar arguments as in \textbf{Part (i)}, we get that for any $n\geq0$,
\begin{align*}
    B_5^{(20t+9)}(4n+2)&\equiv0\pmod2,\\
    B_5^{(20t+9)}(5n+\alpha')&\equiv0\pmod5,
\end{align*}
where $\alpha'\in\{3,4\}$. Combining these, we get the desired result.

\textbf{{Part (iii)}:} Let us take $\ell=5$ and $r=20t+17$, for $t\geq0$, in (\ref{eq:2}) to get
\begin{align*}
    \sum_{n=0}^\infty B_5^{(20t+17)}(n)q^n=\frac{E_5^{20t+17}}{E_1^{20t+17}}.
\end{align*}
The above equation can be rewritten as
\begin{align*}
    \sum_{n=0}^\infty B_5^{(20t+17)}(n)q^n=\frac{E_5^{20t+16}}{E_1^{20t+16}}\cdot\frac{E_5}{E_1}
\end{align*}
as well as
\begin{align*}
    \sum_{n=0}^\infty B_5^{(20t+17)}(n)q^n=\frac{E_5^{20t+17}}{E_1^{20t+20}}\cdot E_1^3.
\end{align*}
Following similar arguments as in \textbf{Part (i)}, we get that for any $n\geq0$,
\begin{align*}
    B_5^{(20t+17)}(4n+2)&\equiv0\pmod2,\\
    B_5^{(20t+17)}(5n+\alpha')&\equiv0\pmod5,
\end{align*}
where $\alpha'\in\{2,3,4\}$. Combining these, we get the desired result.
\end{proof}

\begin{proof}[\textbf{Proof of Theorem \ref{thm5}}]
\textbf{{Part (i)}:} Let us take $\ell=35$ and $r=35t+4$, for $t\geq0$, in (\ref{eq:2}) to get
\begin{align} \label{eq:38}
    \sum_{n=0}^\infty B_{35}^{(35t+4)}(n)q^n=\frac{E_{35}^{35t+4}}{E_{1}^{35t+4}}.
\end{align}
The above equation can be rewritten as
\begin{align*}
    \sum_{n=0}^\infty B_{35}^{(35t+4)}(n)q^n=\frac{E_{35}^{35t+4}}{E_{1}^{35t+5}}\cdot E_1.
\end{align*}
Now, by (\ref{eq:35}), we get
\begin{align*}
    \sum_{n=0}^\infty B_{35}^{(35t+4)}(n)q^n=\frac{E_{35}^{35t+4}}{E_{1}^{35t+5}}\cdot E_{25}\left(\frac{1}{R(q^5)}-q-q^2R\left(q^5\right)\right).
\end{align*}
In view of (\ref{eq:3}), we get
\begin{align*}
    \sum_{n=0}^\infty B_{35}^{(35t+4)}(n)q^n\equiv\frac{E_{35}^{35t+4}E_{25}}{E_{5}^{7t+1}}\left(\frac{1}{R(q^5)}-q-q^2R\left(q^5\right)\right)\pmod5.
\end{align*}
Comparing coefficients of $q^{5n+\alpha_1}$ on both sides, where $\alpha_1\in\{3,4\}$, we get
\begin{align} \label{eq:39}
    B_{35}^{(35t+4)}(5n+\alpha_1)\equiv0\pmod5 \quad \forall \quad n\geq0.
\end{align}

On the other hand, (\ref{eq:38}) can be also rewritten as
\begin{align*}
    \sum_{n=0}^\infty B_{35}^{(35t+4)}(n)q^n=\frac{E_{35}^{35t+4}}{E_{1}^{35t+7}}\cdot E_1^3.
\end{align*}
Now, by (\ref{eq:25}), we get
\begin{align*}
    \sum_{n=0}^\infty B_{35}^{(35t+4)}(n)q^n=\frac{E_{35}^{35t+4}}{E_{1}^{35t+7}}\cdot E_{49}^3\left(\frac{A_2(q^7)}{A_1(q^7)}-q\frac{A_3(q^7)}{A_2(q^7)}-q^2+q^5\frac{A_1(q^7)}{A_3(q^7)}\right)^3.
\end{align*}
In view of (\ref{eq:3}), we get
\begin{align*}
    \sum_{n=0}^\infty B_{35}^{(35t+4)}(n)q^n\equiv\frac{E_{35}^{35t+4}E_{49}^3}{E_{7}^{5t+1}}\Bigg(&\frac{A_2^3}{A_1^3}+4q\frac{A_2A_3}{A_1^2}+4q^2\frac{A_2^2}{A_1^2}+3q^2\frac{A_3^2}{A_1A_2}+6q^3\frac{A_3}{A_1}+6q^3\frac{A_3^3}{A_2^3}\\
    &+3q^4\frac{A_2}{A_1}+4q^4\frac{A_3^2}{A_2^2}+3q^5\frac{A_2^2}{A_1A_3}+4q^5\frac{A_3}{A_2}+q^7\frac{A_2}{A_3}\\
    &+3q^7\frac{A_1A_3}{A_2^2}+6q^8\frac{A_1}{A_2}+3q^9\frac{A_1}{A_3}+3q^{10}\frac{A_1A_2}{A_3^2}+4q^{11}\frac{A_1^2}{A_2A_3}\\
    &+4q^{12}\frac{A_1^2}{A_3^2}+q^{15}\frac{A_1^3}{A_3^3}\Bigg)\pmod7,
\end{align*}
where $A_i:=A_i(q^7)$ for $i\in\{1,2,3\}$. Comparing coefficients of $q^{7n+6}$ on both sides, we get
\begin{align} \label{eq:40}
    B_{35}^{(35t+4)}(7n+6)\equiv0\pmod7 \quad \forall \quad n\geq0.
\end{align}
Combining (\ref{eq:39}) and (\ref{eq:40}), we get the desired result.

\textbf{{Part (ii)}:} Let us take $\ell=35$ and $r=35t+6$, for $t\geq0$, in (\ref{eq:2}) to get
\begin{align*}
    \sum_{n=0}^\infty B_{35}^{(35t+6)}(n)q^n=\frac{E_{35}^{35t+6}}{E_{1}^{35t+6}}.
\end{align*}
The above equation can be rewritten as
\begin{align*}
    \sum_{n=0}^\infty B_{35}^{(35t+6)}(n)q^n=\frac{E_{35}^{35t+6}}{E_{1}^{35t+10}}\cdot E_1^4
\end{align*}
as well as
\begin{align*}
    \sum_{n=0}^\infty B_{35}^{(35t+6)}(n)q^n=\frac{E_{35}^{35t+6}}{E_{1}^{35t+7}}\cdot E_1.
\end{align*}
Following similar arguments as in \textbf{Part (i)}, we get that for any $n\geq0$,
\begin{align*}
    B_{35}^{(35t+6)}(5n+4)&\equiv0\pmod5,\\
    B_{35}^{(35t+6)}(7n+\alpha_2)&\equiv0\pmod7,
\end{align*}
where $\alpha_2\in\{3,4,6\}$. Combining these, we get the desired result.

\textbf{{Part (iii)}:} Let us take $\ell=35$ and $r=35t+11$, for $t\geq0$, in (\ref{eq:2}) to get
\begin{align*}
    \sum_{n=0}^\infty B_{35}^{(35t+11)}(n)q^n=\frac{E_{35}^{35t+11}}{E_{1}^{35t+11}}.
\end{align*}
The above equation can be rewritten as
\begin{align*}
    \sum_{n=0}^\infty B_{35}^{(35t+11)}(n)q^n=\frac{E_{35}^{35t+11}}{E_{1}^{35t+15}}\cdot E_1^4
\end{align*}
as well as
\begin{align*}
    \sum_{n=0}^\infty B_{35}^{(35t+11)}(n)q^n=\frac{E_{35}^{35t+11}}{E_{1}^{35t+14}}\cdot E_1^3.
\end{align*}
Following similar arguments as in \textbf{Part (i)}, we get that for any $n\geq0$,
\begin{align*}
    B_{35}^{(35t+11)}(5n+4)&\equiv0\pmod5,\\
    B_{35}^{(35t+11)}(7n+6)&\equiv0\pmod7.
\end{align*}
Combining these, we get the desired result.

\textbf{{Part (iv)}:} Let us take $\ell=35$ and $r=35t+27$, for $t\geq0$, in (\ref{eq:2}) to get
\begin{align*}
    \sum_{n=0}^\infty B_{35}^{(35t+27)}(n)q^n=\frac{E_{35}^{35t+27}}{E_{1}^{35t+27}}.
\end{align*}
The above equation can be rewritten as
\begin{align*}
    \sum_{n=0}^\infty B_{35}^{(35t+27)}(n)q^n=\frac{E_{35}^{35t+27}}{E_{1}^{35t+30}}\cdot E_1^3
\end{align*}
as well as
\begin{align*}
    \sum_{n=0}^\infty B_{35}^{(35t+27)}(n)q^n=\frac{E_{35}^{35t+27}}{E_{1}^{35t+28}}\cdot E_1.
\end{align*}
Following similar arguments as in \textbf{Part (i)}, we get that for any $n\geq0$,
\begin{align*}
    B_{35}^{(35t+27)}(5n+\alpha_1)&\equiv0\pmod5,\\
    B_{35}^{(35t+27)}(7n+\alpha_2)&\equiv0\pmod7,
\end{align*}
where $\alpha_1\in\{2,3,4\}$ and $\alpha_2\in\{3,4,6\}$. Combining these, we get the desired result.

\textbf{{Part (v)}:} Let us take $\ell=35$ and $r=35t+32$, for $t\geq0$, in (\ref{eq:2}) to get
\begin{align*}
    \sum_{n=0}^\infty B_{35}^{(35t+32)}(n)q^n=\frac{E_{35}^{35t+32}}{E_{1}^{35t+32}}.
\end{align*}
The above equation can be rewritten as
\begin{align*}
    \sum_{n=0}^\infty B_{35}^{(35t+32)}(n)q^n=\frac{E_{35}^{35t+32}}{E_{1}^{35t+35}}\cdot E_1^3
\end{align*}
Following similar arguments as in \textbf{Part (i)}, we get that for any $n\geq0$,
\begin{align*}
    B_{35}^{(35t+32)}(5n+\alpha_1)&\equiv0\pmod5,\\
    B_{35}^{(35t+32)}(7n+6)&\equiv0\pmod7,
\end{align*}
where $\alpha_1\in\{2,3,4\}$. Combining these, we get the desired result.

\textbf{{Part (vi)}:} Let us take $\ell=35$ and $r=35t+34$, for $t\geq0$, in (\ref{eq:2}) to get
\begin{align*}
    \sum_{n=0}^\infty B_{35}^{(35t+34)}(n)q^n=\frac{E_{35}^{35t+34}}{E_{1}^{35t+34}}.
\end{align*}
The above equation can be rewritten as
\begin{align*}
    \sum_{n=0}^\infty B_{35}^{(35t+34)}(n)q^n=\frac{E_{35}^{35t+34}}{E_{1}^{35t+35}}\cdot E_1
\end{align*}
Following similar arguments as in \textbf{Part (i)}, we get that for any $n\geq0$,
\begin{align*}
    B_{35}^{(35t+34)}(5n+\alpha_1)&\equiv0\pmod5,\\
    B_{35}^{(35t+34)}(7n+\alpha_2)&\equiv0\pmod7,
\end{align*}
where $\alpha_1\in\{3,4\}$ and $\alpha_2\in\{3,4,6\}$. Combining these, we get the desired result.
\end{proof}

\begin{proof}[\textbf{Proof of Theorem \ref{thm6}}]
\textbf{{Part (i)}:} Let us take $\ell=55$ and $r=55t+21$, for $t\geq0$, in (\ref{eq:2}) to get
\begin{align} \label{eq:41}
    \sum_{n=0}^\infty B_{55}^{(55t+21)}(n)q^n=\frac{E_{55}^{55t+21}}{E_{1}^{55t+21}}.
\end{align}
The above equation can be rewritten as
\begin{align*}
    \sum_{n=0}^\infty B_{55}^{(55t+21)}(n)q^n=\frac{E_{55}^{55t+21}}{E_{1}^{55t+25}}\cdot E_1^4.
\end{align*}
Now, by (\ref{eq:35}), we get
\begin{align*}
    \sum_{n=0}^\infty B_{55}^{(55t+21)}(n)q^n=\frac{E_{55}^{55t+21}}{E_{1}^{55t+25}}\cdot E_{25}^4\left(\frac{1}{R(q^5)}-q-q^2R\left(q^5\right)\right)^4.
\end{align*}
In view of (\ref{eq:3}), we get
\begin{align*}
    \sum_{n=0}^\infty B_{55}^{(55t+21)}(n)q^n\equiv\frac{E_{55}^{55t+21}E_{25}^4}{E_{5}^{11t+5}}\Bigg(&\frac{1}{R\left(q^5\right)^4}+\frac{q}{R\left(q^5\right)^3}+\frac{2q^2}{R\left(q^5\right)^2}+\frac{3q^3}{R\left(q^5\right)}+2q^5R\left(q^5\right)\\
    &+2q^6R\left(q^5\right)^2+4q^7R\left(q^5\right)^3+q^8R\left(q^5\right)^4\Bigg)\pmod5.
\end{align*}
Comparing coefficients of $q^{5n+4}$ on both sides, we get
\begin{align} \label{eq:43}
    B_{55}^{(55t+21)}(5n+4)\equiv0\pmod5 \quad \forall \quad n\geq0.
\end{align}

On the other hand, (\ref{eq:41}) can be also rewritten as
\begin{align*}
    \sum_{n=0}^\infty B_{55}^{(55t+21)}(n)q^n=\frac{E_{55}^{55t+21}}{E_{1}^{55t+22}}\cdot E_1.
\end{align*}
Now, by (\ref{eq:42}), we get
\begin{align*}
    \sum_{n=0}^\infty B_{55}^{(55t+21)}(n)q^n=\frac{E_{55}^{55t+21}}{E_{1}^{55t+22}}\cdot E_{121}\Bigg(&\frac{B_4(q^{11})}{B_2(q^{11})}-q\frac{B_2(q^{11})}{B_1(q^{11})}-q^2\frac{B_5(q^{11})}{B_3(q^{11})}+q^5+q^7\frac{B_3(q^{11})}{B_4(q^{11})}\\
    &-q^{15}\frac{B_1(q^{11})}{B_5(q^{11})}\Bigg).
\end{align*}
In view of (\ref{eq:3}), we get
\begin{align*}
    \sum_{n=0}^\infty B_{55}^{(55t+21)}(n)q^n\equiv\frac{E_{55}^{55t+21}E_{121}}{E_{11}^{5t+2}}\Bigg(&\frac{B_4(q^{11})}{B_2(q^{11})}-q\frac{B_2(q^{11})}{B_1(q^{11})}-q^2\frac{B_5(q^{11})}{B_3(q^{11})}+q^5+q^7\frac{B_3(q^{11})}{B_4(q^{11})}\\
    &-q^{15}\frac{B_1(q^{11})}{B_5(q^{11})}\Bigg)\pmod{11}.
\end{align*}
Comparing coefficients of $q^{11n+\alpha_2}$ on both sides, where $\alpha_2\in\{3,6,8,9,10\}$, we get
\begin{align} \label{eq:44}
    B_{55}^{(55t+21)}(11n+\alpha_2)\equiv0\pmod{11} \quad \forall \quad n\geq0.
\end{align}
Combining (\ref{eq:43}) and (\ref{eq:44}), we get the desired result.

\textbf{{Part (ii)}:} Let us take $\ell=55$ and $r=55t+32$, for $t\geq0$, in (\ref{eq:2}) to get
\begin{align*}
    \sum_{n=0}^\infty B_{55}^{(55t+32)}(n)q^n=\frac{E_{55}^{55t+32}}{E_{1}^{55t+32}}.
\end{align*}
The above equation can be rewritten as
\begin{align*}
    \sum_{n=0}^\infty B_{55}^{(55t+32)}(n)q^n=\frac{E_{55}^{55t+32}}{E_{1}^{55t+35}}\cdot E_1^3
\end{align*}
as well as
\begin{align*}
    \sum_{n=0}^\infty B_{55}^{(55t+32)}(n)q^n=\frac{E_{55}^{55t+32}}{E_{1}^{55t+33}}\cdot E_1.
\end{align*}
Following similar arguments as in \textbf{Part (i)}, we get that for any $n\geq0$,
\begin{align*}
    B_{55}^{(55t+32)}(5n+\alpha_1)&\equiv0\pmod5,\\
    B_{55}^{(55t+32)}(11n+\alpha_2)&\equiv0\pmod{11},
\end{align*}
where $\alpha_1\in\{2,3,4\}$ and $\alpha_2\in\{3,6,8,9,10\}$. Combining these, we get the desired result.

\textbf{{Part (iii)}:} Let us take $\ell=55$ and $r=55t+54$, for $t\geq0$, in (\ref{eq:2}) to get
\begin{align*}
    \sum_{n=0}^\infty B_{55}^{(55t+54)}(n)q^n=\frac{E_{55}^{55t+54}}{E_{1}^{55t+54}}.
\end{align*}
The above equation can be rewritten as
\begin{align*}
    \sum_{n=0}^\infty B_{55}^{(55t+54)}(n)q^n=\frac{E_{55}^{55t+54}}{E_{1}^{55t+55}}\cdot E_1.
\end{align*}
Following similar arguments as in \textbf{Part (i)}, we get that for any $n\geq0$,
\begin{align*}
    B_{55}^{(55t+54)}(5n+\alpha_1)&\equiv0\pmod5,\\
    B_{55}^{(55t+54)}(11n+\alpha_2)&\equiv0\pmod{11},
\end{align*}
where $\alpha_1\in\{3,4\}$ and $\alpha_2\in\{3,6,8,9,10\}$. Combining these, we get the desired result.
\end{proof}
\bibliographystyle{abbrv}

\bibliography{bibliography}

\end{document}